\documentclass[11pt]{article}
                \usepackage{amsmath, amsthm, graphics, amssymb,fullpage, color,setspace, graphicx}

                \newtheorem* {Main} {Theorem 1}
                
                \newtheorem* {Prop1} {Proposition 5.1 (the upper bound)}
                \newtheorem* {Prop2} {Proposition 5.2 (the lower bound)}
                
                \newtheorem* {twoactions} {Lemma 5.5}
                
                \newtheorem* {wordgrowth}{Proposition 2.3}
                \newtheorem* {lengthLemma} {Lemma 2.2}
                
                \newtheorem* {representation} {Lemma 5.6}
                \newtheorem* {bruhattits} {Proposition 5.3}

                \newtheorem* {algebraicentropy} {Proposition 3.1}
                \newtheorem* {lowerbound} {Lemma 5.4}
                \newtheorem* {specrad} {Proposition 2.1}
                
\begin{document}

                \title{Algebraic Entropy and the Action of Mapping Class Groups on Character Varieties}
                \author{Asaf Hadari}
                \maketitle

                \begin{abstract}
                We extend the definition of algebraic entropy to endomorphisms of affine varieties. We then calculate the algebraic entropy of the action of elements of mapping class groups on various character varieties, and show that it is equal to a quantity we call the spectral radius, a generalization of the dilatation of a Pseudo-Anosov mapping class. Our calculations are compatible with all known calculations of the topological entropy of this action.
                \end{abstract}

                \section{Introduction}

                Let $S = S_{g,b}$ be an oriented surface of genus $g$ with $b \geq 1$ boundary components. The $\textit{Mapping class group of S}$, which we denote by $\textrm{Mod}(S)$, consists of isotopy classes of orientation-preserving diffeomorphisms of $S$ which fix the boundary components pointwise.

                Choose a basepoint $p_{0} \in S$, and let $\pi = \pi_{1}(S,p_{0})$. Given an algebraic group $G$, one can construct the variety:                 $Hom(\pi,G)$, which is called  the \textit{G-representation variety of} $\pi$. The group $G$ acts algebraically on this variety by conjugation. The categorical quotient of the representation variety is called the \textit{G-character variety of} $\pi$, which we will often denote by $\mathfrak{X} := Hom(\pi, G)//G$.

                The mapping class group $\textrm{Mod}(S)$ acts algebraically on $\mathfrak{X}$. Given $f \in \textrm{Mod}(S)$, our goal in this paper is to calculate an algebraic invariant that gives a measure of the complexity of the action of $f$ on $\mathfrak{X}$.

                In the study of dynamical systems, there are several different measures of complexity called entropy - topological entropy and measure theoretic entropy being two common examples. In general, one expects actions with high entropy to be more complicated than actions with low entropy. For our purposes we wish to use an measure that captures the algebraic nature of the action. In \cite{belvial}, Bellon and Viallet  define a notion called algebraic entropy for algebraic endomorphisms of affine space, which measures the growth rate of the degrees of iterates of the map.

                The variety $\mathfrak{X}$ is affine, but there is no preferred  way to embed it into affine space. One of the goals of this paper is to give an intrinsic natural extension of Bellon and Viallet's concept of algebraic entropy to algebraic self maps of affine varieties. This is the invariant we study. As a caution to the reader, we mention that there is a different dynamical invariant, due to Gromov, which is called algebraic entropy. We define all the terms we use, so no confusion should arise.

                Let $\textrm{e}_{\textrm{alg}}(f)$ be the algebraic entropy of $f$, which is defined in
                Section 3, and let $\rho(f)$ is the \textit{spectral radius of f}, a
                generalization of the log of the dilatation of a pseudo-Anosov element, which
                is defined in Section 2. We prove the following theorem.

                \begin{Main}
                Let $K = \mathbb{R}$ or $\mathbb{C}$ and $G$ be one of the following groups:
                $$SL_{N}(K), GL_{N}(K), O_{N}(\mathbb{R}) (N \geq 3), SO_{N}(\mathbb{R}) (n \geq 3), U_{N}, SU_{2}, Sp_{2N}(\mathbb{R})$$

                Let $S$ be a surface with free fundamental group, and let $f \in \textrm{Mod}(S)$. The mapping class $f$ acts on the $G$ character variety of $S$, and one has that $$\textrm{e}_{\textrm{alg}}(f) = \rho(f)$$
                \end{Main}

                The topological entropy of mapping class group actions on character
                varieties has been calculated by Fried for the case $S = S_{1,1}$
                and $G = SU(2)$ \cite{entcharfried} and by Cantat and Loray for reduced character varieties (these are character varieties where the traces of boundary components are fixed) in the case $S =
                S_{0,4}$, $G = SL_{2}(\mathbb{C})$ \cite{CantatLoray}. The algebraic
                entropy was calculated by Brown for the case $S = S_{1,1}$ and $G =
                SU(2)$ and a specific embedding of $\mathfrak{X}$ \cite{brownentropy}. In
                all of the above cases, the entropy calculated was equal to
                $\rho(f)$.

                The paper is organized as follows. In section $2$ we define the concept of spectral radius and show how to calculate it for many elements of the mapping class group. In section $3$ we define the concept of algebraic entropy. In section $4$ we discuss the basics of character varieties and define the action of mapping class groups on them. Section $5$ is devoted to the proof of theorem $1$, divided into the proof of two inequalities.

                \paragraph{Acknowledgements.} The author wishes to thank Khalid Bou-Rabee, Thomas Zamojski, and Benson Farb for their comments and for many illuminating discussions. He also wishes to thank the anonymous referee for extensive and invaluable comments.

                \section{The spectral radius of a mapping class}

                \paragraph{Mapping class groups.} Let $S = S_{g,b}$ be a surface of genus $g$ with $b$ boundary
                components (in this paper we will always assume that $b \geq 1$). Let $\textrm{Diff}^{+}(S)$ be the group of orientation-preserving diffeomorphisms of $S$ that are the identity on the
                boundary components. The \textit{mapping class group of S} is the
                group $\textrm{Mod}(S) = \pi_{0}(\textrm{Diff}^{+}(S))$.

                In what follows we will make no notational distinctions between simple closed curves and their homotopy classes. Also, we will assume that a base point is chosen on the boundary of $S$. This allows us to identify any $f \in \textrm{Mod}(S)$ with an element of $\textrm{Aut}(\pi)$.  All of the information about mapping class groups used in this paper can be found in \cite{mcg}.

                \paragraph{Definition} Let $f \in \textrm{Mod}(S)$. Let $\mathcal{S}$ be a generating set for $\pi$. For $w \in \pi$,
                let $|w|_{\mathcal{S}, red}$ be its cyclically reduced word length with respect to $\mathcal{S}$. Define the \textit{spectral radius of f with respect to} $\mathcal{S}$ to be the quantity:
                $$\rho^{\mathcal{S}}(f): = \sup_{\alpha \in \pi} \textit{ } \limsup_{n \to \infty} \frac{1}{n} \log(|f^{\circ n} \alpha|_{\mathcal{S},red})$$

                First notice that the above definition does not depend on the choice of the base point. Indeed, after changing the base point, the action of $f$ on $\pi$ is changed by composition with an inner automorphism. This clearly does not change cyclically reduced word lengths. The next proposition shows that the dependence on the set $\mathcal{S}$ can be dropped.

                \begin{specrad} Given any element $f \in \textrm{Mod}(S)$, and any two generating sets $\mathcal{S}_{1}$ and $\mathcal{S}_{2}$ of $\pi$, the following equality holds:
                $$\rho^{\mathcal{S}_{1}}(f) = \rho^{\mathcal{S}_{2}}(f) $$
                \end{specrad}

                \paragraph{Proof.} Recall that a map $\Phi: X \rightarrow Y$ between metric spaces is
                called a \textit{quasi-isometry} if there exist positive constants
                $K,C, D$ such that for every $x,y \in X$:
                $$\frac{1}{K}d_{X}(x,y) - C \leq d_{Y}(f(x),f(y)) \leq Kd_{X}(x,y)
                +C
                 $$
                and such that for every $w \in Y$, $d_{Y}(w, \Phi(X)) \leq D$.

                The generating set $\mathcal{S}_{1}$ and $\mathcal{S}_{2}$ define two word metrics on $\pi$. It is well known that the two metric spaces defined in this way are quasi-isometric. Any element $w \in \pi$ acts on $\pi$ by left translation. It is well known that the translation length of this action in the metric given by $|.|_{\mathcal{S}_{i}}$ ($i = 1,2$) is $|w|_{\mathcal{S}_{i},red}$.

                Suppose $K,C,D$ are the quasi-isometry constants for the quasi-isomorphism between $(\pi, |.|_{\mathcal{S}_{1}})$ and $(\pi, |.|_{\mathcal{S}_{2}})$. Using the characterization of $|w|_{\mathcal{S}_{i},red}$ as a translation length, it is clear that:
                $$\frac{1}{K}|w|_{\mathcal{S}_{1},red} + C \leq |w|_{\mathcal{S}_{1},red} \leq K|w|_{\mathcal{S}_{1},red} + C $$

                Given any constants $A,B$ it's true that $$\sup_{\alpha \in \pi} \textit{ } \limsup_{n \to \infty} \frac{1}{n} \log(A|f^{\circ n} \alpha|_{\mathcal{S}_{i},red} + B) = \sup_{\alpha \in \pi} \textit{ } \limsup_{n \to \infty} \frac{1}{n} \log(|f^{\circ n} \alpha|_{\mathcal{S}_{i},red})$$

                And thus $\rho^{\mathcal{S}_{1}}(f) = \rho^{\mathcal{S}_{2}}(f) $, as required. $\Box$

               \bigskip

                Since the definition of spectral radius does not depend on the generating set, we will suppress the $\mathcal{S}$ in the notation, and use $\rho(f)$ for the \textit{spectral radius of} $f$.

                \paragraph{Calculating the spectral radius.} Our next goal is to calculate spectral radius for many elements of the mapping class group.

                A \textit{multicurve} in $S$ is a finite collection of homotopy classes of mutually disjoint simple closed curves in $S$, none of which is homotopic to a boundary component. The mapping class group acts on the set of multicurves. This action can be used to classify elements of the mapping class group as follows.  \\
                Let $f \in
                \textrm{Mod}(S)$. Exactly one of the following is true.
                \begin{enumerate}
                \item The order of $f$ is finite.
                \item The order of $f$ is infinite, and there exists a multicurve $M$ such that $f(M) = M$. In this case $f$ is called \textit{reducible}.
                \item For every multicurve,$M$, in $S$, $f(M) \neq M$. In this case $f$ is
                called \textit{pseudo-Anosov}.
                \end{enumerate}

                We will consider a particularly well behaved subclass of $\textrm{Mod}(S)$, called pure elements. We say that an element $f \in \textrm{Mod}(S)$ is \textit{pure} if there exists a diffeomorphism $\phi$ of $S$ in the homotopy class of $f$, and a (possibly empty) one dimensional submanifold $c$ of $S$ with the following properties.

                \begin{enumerate}
                \item None of the components of $c$ are null-homotopic or homotopic to boundary components of $S$.
                \item $\phi|_{c} = id$
                \item $\phi$ does not rearrange the components of $S \backslash c$.
                \item On each component of $S_{c}$, the surface obtained by cutting $S$ along $c$, $\phi$ induces a diffeomorphism which is homotopic either to the identity or to a pseudo-Anosov.
                \end{enumerate}

                Note that any pseudo-Anosov mapping class is pure. In \cite{crs}, it is proved that $\textrm{Mod}(S)$ contains a finite index subgroup consisting entirely of pure elements. An example of such a group is the kernel of the action of $\textrm{Mod}(S)$ on $H_{1}(S, \mathbb{Z} / 3\mathbb{Z})$.

                A theorem of Thurston describes a canonical geometric element contained in a
                pseudo-Anosov mapping class $f$. Using this element, one can attach
                an algebraic integer $\lambda = \lambda(f) > 1$ to $f$ called the
                \textit{dilatation of f}. To any pure element, one can attach a collection of dilatations, one for each component of $S_{c}$ on which $f$ acts as a pseudo-Anosov. We call the maximum one of these dilatations the \textit{dilatation of} $f$, and denote it $\lambda(f)$. We now proceed to calculate the spectral radius of any pure element of the mapping class group.

  \begin{lengthLemma}
                Let $f \in \textrm{Mod}(S)$ be a pure element, and let $\alpha$ be the isotopy class of a simple
                closed curve on $S$. Let $g$ be a Riemannian metric on $S$. If we
                denote by $l_{g}(.)$ the $g$-length of a curve of an isotopy class of curves. Then:
                $$\limsup_{n \to \infty} \frac{1}{n} \log(l_{g}(f^{\circ n}(\alpha))) \leq \log(\lambda(f)) $$
                Furthermore, there exists a simple closed curve $\alpha$ for which
                the above inequality is an equality.
                \end{lengthLemma}

                \paragraph{Proof.}
                When $f$ is a pseudo-Anosov element with dilatation $\lambda$, then
                $\rho(f) = \log(\lambda)$, and the inequality in the claim of the
                lemma is an equality for every curve. The proof of this fact can be
                found for instance in \cite{mcg} (Theorem $13.20$).

                Suppose that $f$ is reducible. Let $c$ and $\phi$ be the one dimensional submanifold, and the diffeormorphism associated with $f$. Suppose that $\alpha$ is transverse to $c$. The finite set $\alpha \cap c$ is fixed by $\phi$. Thus, $\alpha \backslash \phi$ consists of a collection of arcs: $\alpha_{1} \ldots \alpha_{p}$ whose endpoints are fixed by $\phi$. We view each of these arcs as being a subset of one of the components of the surface obtained from $S$ by cutting along $c$. The result quoted in \cite{mcg} can be restated to apply to arcs whose endpoints are fixed by $f$. The same proof carries through with only notational changes. One has that:

                $$l_{g}(\phi^{\circ n}(\alpha)) \leq \sum_{i = 1}^{p}l_{g}(\phi^{\circ n}(\alpha_{i}))  $$

                The first part of the result is now clear. To see the second part, choose a curve contained in a subsurface with boundary on which the action of $\phi$ has dilatation $\lambda(f)$.  $\Box$

     \begin{wordgrowth}
                Let $f \in \textrm{Mod}(S)$ be a pure element. Then $$\rho(f) = \log(\lambda(f)) $$
                \end{wordgrowth}

                \paragraph{Proof.}

                Choose a generating set $\mathcal{S}$ of $\pi$. This  defines a word metric on $\pi$. A choice of a hyperbolic metric $g$ on $S$, and
                a  basepoint $x_{0} \in \mathbb{H}^{2}$, the upper half plane, defines an embedding of the Cayley graph into the hyperbolic plane. This embedding
                induces a new metric on the graph. It is well known that these two
                metrics are quasi-isometric. Every element $\alpha \in \pi$ acts as
                an isometry on the Caley graph of $\pi$. Furthermore, $\alpha$ acts
                as an isometry on $\mathbb{H}^{2}$ which preserves the embedded
                Cayley graph. The actions of $\alpha$ on the embedded Cayley graph
                and the abstract Cayley graph are conjugate.

                For every $n$, and every curve $\gamma$, the curve $f^{\circ n}(\gamma)$ corresponds to an isometry
                on $\mathbb{H}^{2}$. The translation length of this isometry is
                given by $l_{g}(f^{\circ n}(\gamma))$. The curve $f^{\circ n}(\gamma)$ also acts on the
                Cayley graph of $\pi$ by left multiplication. The translation length
                of that action is $|f^{\circ n}(\gamma)|_{\mathcal{S}, red}$. The first
                assertion of the lemma now clearly follows from quasi-isometry, and
                the fact that:
                $$\lim_{n \to \infty} \frac{1}{n}
                \log(l_{g}(f^{n\circ }(\gamma))) \leq \rho(f)$$
                The second assertion
                follows immediately from the second part of Lemma 2.2.  $\Box$

                \section{Algebraic entropy}

                Let $f$ be an endomorphism of $\mathbb{A}^{N}$. Define the
                \textit{dynamical degree of f} as:
                $$\Delta(f) = \lim_{n \to \infty} \deg(f^{\circ n})^{\frac{1}{n}} $$
                In \cite{belvial} Bellon and Viallet define the \textit{algebraic entropy} of $f$ as $$\textrm{e}_{\textrm{alg}}(f) =
                \log(\Delta(f))$$

                Algebraic entropy is meant to be an algebraic approximation of topological entropy. To see this, consider the following heuristic argument: the topological entropy of $f$ can often be estimated by
                calculating the exponential growth rate of the number of isolated
                fixed points of $f^{\circ n}$. For a polynomial automorphism $f$ of
                $\mathbb{C}^{m}$, the number of isolated fixed points of $f$ is at
                most $\deg(f)$. Thus, calculating the exponential growth rate of the
                degrees of $f^{n}$ can be seen as estimating the topological entropy
                of $f$.

                In this paper we are concerned with endomorphisms of character varieties, which are affine varieties. Given a variety $V$ equipped with an endomorphism $f$, we wish to give a definition that provides an algebraic approximation of the topological entropy of $f$. The naive approach is to embed $V$ to affine space, extend $f$ to an endomorphism of affine space and calculate its algebraic entropy. The problem with this approach is that neither the embeddings nor the extensions are canonical, and one can get many different results in this way. Our goal in this section is to give an intrinsic invariant which generalizes algebraic entropy.

                \paragraph{Definition} Let $V \subset \mathbb{A}^{N}$ be a subvariety of affine $N$-space. Let $f : \mathbb{A}^{N} \to \mathbb{A}^{M}$ be a morphism. Let $\mathfrak{R}_{f,V}$ be the set of morphisms, $g: \mathbb{A}^{N} \to \mathbb{A}^{M}$ such that $g|_{V} = f|_{V}$.  Define the \textit{degree of f relative to V} as the quantity:
                $$\deg(f;V) := \min_{g \in \mathfrak{R}_{f,V}}\deg(g) $$

                \paragraph{Definition} Let $V \subset \mathbb{A}^{N}$ be a subvariety of affine $N$-space. Let $f : \mathbb{A}^{N} \to \mathbb{A}^{N}$ be a morphism. Define the \textit{algebraic entropy of f relative to V} as the quantity:
                $$\textrm{e}_{\textrm{alg}}(f;V) =  \limsup_{n \to \infty} \frac{1}{n} \log \deg(f^{\circ n};V)  $$

                Suppose now that $V$ is an affine variety, and $f: V \to V$ is a morphism. $V$ can be embedded in many ways into affine space, and $f$ can be extended in many ways to a morphism of affine space. For each embedding of $V$ and each extension of $f$ we can calculate the algebraic entropy relative to $V$. The next proposition shows that the above choices don't affect the algebraic entropy.

                \begin{algebraicentropy}
                Let $V$ be an affine variety, and let $f: V \to V$ be a morphism. Let $\iota_{1}: V \to \mathbb{A}^{N_{1}}$ and $\iota_{2}: V \to \mathbb{A}^{N_{2}}$ be two affine embeddings of $V$. Let $g_{i}: \mathbb{A}^{N_{i}} \to \mathbb{A}^{N_{i}}$ ($i = 1,2$) be morphisms such that: $g_{i}(\iota_{i}(V)) = \iota_{i}(V) $, and $\iota^{*}_{i}(g_{i}) = f$. Then:
                $$\textrm{e}_{\textrm{alg}}(g_{1};\iota_{1}(V)) = \textrm{e}_{\textrm{alg}}(g_{2};\iota_{2}(V)) $$
                \end{algebraicentropy}

            \paragraph{Proof.} First note that if $\iota_{1} = \iota_{2}$ then the claim is trivial by the definition of relative algebraic entropy.
\\

            For $i = 1,2$ the maps $\iota_{i}$ can be written in coordinates as:
            $$\iota_{i} = (x_{j,i})_{j = 1}^{N_{i}}$$
             Since each $\iota_{i}$ is an embedding, then each of the sets $X_{i} = \{x_{1,i}, \ldots, x_{N_{i},i} \}$ generates the ring of functions of $V$. Thus, for each $j = 1, \ldots, N_{2}$, we can non-canonically write $x_{j,2}$ as a polynomial in the elements of $X_{1}$. Using this, we get a morphism $p_{1}: \mathbb{A}^{N_{1}} \to \mathbb{A}^{N_{2}}$, such that $p_{1} \circ \iota_{1} = \iota_{2}$. Say that the degree of $p_{1}$ is $D_{1}$. By similar reasoning, there is a morphism $p_{2}: \mathbb{A}^{N_{2}} \to \mathbb{A}^{N_{1}}$ of degree $D_{2}$ such that $p_{1} \circ \iota_{2} = \iota_{1}$.
             \smallskip
             Given an endomorphism $\tau: \mathbb{A}^{N_{2}} \to \mathbb{A}^{N_{2}}$ of degree $t$, we can construct the endomorphism $p_{2} \circ \tau \circ p_{1}: \mathbb{A}^{N_{1}} \to \mathbb{A}^{N_{1}}$. The degree of this morphism is at most $tD_{1}D_{2}$.
             \\

             For any integer $n$, we take $\tau = g_{2}^{\circ n}$. The resulting endomorphism is clearly an extension of $f^{\circ n}$ from $\iota_{1}(V)$ to $\mathbb{A}^{N_{1}}$. By the definition of degree relative to a subvariety, we have that:
            $$\deg(g_{1}^{\circ n}; \iota_{1}(V)) \leq D_{1}D_{2} (\deg(g_{2}^{\circ n}; \iota_{2}(V))) $$
            Thus
            $$\textrm{e}_{\textrm{alg}}(g_{1};\iota_{1}(V)) \leq \textrm{e}_{\textrm{alg}}(g_{2};\iota_{2}(V)) $$
            Reversing the roles played by the two spaces, we get the result. $\Box$

            Using Proposition 3.1, we can now define an intrinsic notion of algebraic entropy.
            \paragraph{Definition} Let $V$ be an affine variety, and let $f: V \to V$ be a morphism. Let $\iota: V \to \mathbb{A}^{N}$ be an affine embedding and let $g: \mathbb{A}^{N} \to \mathbb{A}^{N}$ be a morphism such that $g(\iota(V)) = \iota(V)$ and $\iota^{*}(g) = f$. Define the \textit{algebraic entropy of f} to be the quantity:
            $$\textrm{e}_{\textrm{alg}}(f) = \textrm{e}_{\textrm{alg}}(g;\iota(V)) $$

            Notice that if $V = \mathbb{A}^{N}$, then the above definition agrees with the regular definition of algebraic entropy.

\section{The mapping class group action on character varieties}

                \paragraph{Representation varieties and character varieties.} Suppose that $\pi \cong F_{n}$. Let $G$ be a linear reductive algebraic group defined over the field $K$. Fix, once and for all, a faithful linear representation of $G$.
                Let
                $$\mathfrak{R} = \mathfrak{R}(S,G) = Hom(\pi, G) \cong
                G^{n}$$ The set $\mathfrak{R}$ has a natural structure as a variety. We
                call $\mathfrak{R}$ the $G$ \textit{representation variety of} $S$.

                $G$ acts algebraically on $\mathfrak{R}$ by componentwise
                conjugation. Consider the ring of invariants under this action,
                $F[\mathfrak{R}]^{G}$. Define:
                $$\mathfrak{X} = \mathfrak{X}(S,G) = \mathfrak{R}//G := \textrm{spec}(F[\mathfrak{R}]^{G})$$
                We call this variety the $G$ \textit{character variety of} $S$. We
                think of it as the set of characters of representations of $\pi$
                into $G$.

                The group $\textrm{Aut}(\pi)$ acts on $\mathfrak{R}$ in the following way: given a representation $\phi \in \mathfrak{R}$, an element $\alpha \in \pi$ and an automorphism $f \in \textrm{Aut}(\pi)$, define:

                $$f(\phi)(\alpha) := \phi(f(\alpha)) $$

                Let $\textrm{Out}(\pi)$ be the group of outer automorphisms of $\pi$. The action of $\textrm{Aut}(\pi)$ on $\mathfrak{R}$ descends to an action of $\textrm{Out}(\pi)$ on $\mathfrak{X}$. Since $\textrm{Mod}(S)$ can be viewed as a subgroup of $\textrm{Out}(\pi)$, we get an action of $\textrm{Mod}(S)$ on $\pi$.

                \paragraph{Generating the ring of invariants.} A theorem which gives an explicit generating set for the ring $K[\mathfrak{R}]^{G}$ is often called a \textit{first fundamental theorem for G-invariants of n matrices}, where $n$ is the rank of $\pi$. A first fundamental theorem for $SL_{2}(\mathbb{C})$, $SL_{2}(\mathbb{R})$, and $SU_{2}$ is known since the work of Fricke.  In \cite{tracefunc}, Procesi proves a first fundamental theorem of $GL_{N}$, $SL_{N}$,  $O_{N}$, $U(N)$, and  $Sp_{2N}(\mathbb{R})$ for $m$ matrices. In \cite{Son} the Rogara proves a first fundamental theorem of $SO_{N}(\mathbb{R})$ invariants for $n$ matrices.

                The most common functions that are given as generators are called \textit{trace functions}. Given an element $\alpha \in \pi$, we can define a function $\textrm{tr}_{\alpha}$ on $\mathfrak{R}$ by
                $$\textrm{tr}_{\alpha}(\phi) = \textrm{trace}(\phi(\alpha))$$

                Choosing a generating set $\mathcal{A} = \{X_{1}, \ldots X_{n}\}$ for $\pi$ identifies $\mathfrak{R}$ as a subset of $M_{N \times N}^{n} \cong \mathbb{A}^{N^{2}n}$. Under this identification, to any word $w$ in the elements of $\mathcal{A}$ one can associate the function $\textrm{tr}_{w}$, which is a homogeneous polynomial. Note that given an element $\alpha \in \pi$, it may be possible to write $\alpha$ in several different ways as a word in the elements of $\mathcal{A}$, and thus $\textrm{tr}_{\alpha}$ can be extended in more than one way to a function on $M_{N \times N}^{n}$.

                Since the functions $\textrm{tr}_{\alpha}$ are conjugation invariant, we can view
                $\textrm{tr}_{\alpha}$ as an element of $\mathbb{C}[\mathfrak{X}]$, i.e. as a
                regular function on $\mathfrak{X}$. For $GL_{N}$, the set of trace functions generate the ring of invariants.
                In fact,  only finitely many trace functions are required to generate the ring.  For the other cases, slightly more complicated functions are needed. For example, for the case $G = O_{N}$, one needs to take traces of words in the elements of $\mathcal{A}$ and their transposes. For $SP_{2N}$, one needs to add symplectic transposes.  These functions are all  homogeneous functions on the coordinates of the matrices representing elements of $\pi$. Formally, we use the following fact:

                For any generating set $\mathcal{A}$, there exists an integer $L$, finitely many functions: $h_{i}: F_{\mathcal{A}} \to K[X_{1}, \ldots, X_{N^{2}n}] \textit{   } (i = 1 \ldots p)$ whose images are all homogeneous of degree at most $L$, and a finite subset $\{w_{1}, \ldots, w_{p}\} \subset \pi$ such that:
                \begin{enumerate}
                \item For any $i$, and $\alpha \in \pi$, the function $h_{i}(\alpha)$ is homogeneous of degree at most $L$ in the coordinates of the matrix representing $\alpha$, and is invariant under conjugation.
                \item The collection $h_{i}(\alpha_{i})$ generates the ring of invariants.
                \end{enumerate}

                On a first reading, we suggest that the reader think of all of the functions $h_{i}$ as being the trace function, and the $h_{i}(w_{i})$ as being the traces of finitely many words. In this paper we use these functions to find affine embeddings of $\mathfrak{X}$.

                \paragraph{Example.} Let $S = S_{1,1}$ be a one holed torus, and let $G= SL_{2}(\mathbb{R})$. Choose $2$ simple closed curves on $S$ whose intersection number is $1$. Call these curves $X$ and $Y$, and let $Z = XY$. It is well known that the map $Tr: \mathfrak{X} \to \mathbb{R}^{3}$ given by

                $$Tr(\chi) = (\textrm{tr}_{X}(\chi), \textrm{tr}_{Y}(\chi), \textrm{tr}_{Z}(\chi)) $$

                is an isomorphism. The action of the Dehn twist about $X$ on $\pi$ is given by

                $$T_{X}(X) = X \textit{,              }  T_{X}(Y) = YX $$

                In the above trace coordinates, the action is given by:

                $$(x,y,z) \to (x, z, xz-y) $$

                Now consider the action of $T_{X}^{\circ 2}$. In coordinates, we can write it out as:

$$(x,y,z) \to (x, xz - y, x^{2}z - xy - z) $$

                The action of $T_{X}^{\circ 3}$ is given by:
                $$(x,y,z) \to (x, x^{2}z - yx - z, x^{3}z - x^{2}y + y  ) $$

                In general, it is simple to see that $\deg(T_{X}^{\circ n} = n)$, and thus $\textrm{e}_{\textrm{alg}}(T_{X}) = 0$. This agrees with the fact that $\rho(T_{X}) = 0$.

                \section{Proof of Theorem 1.}

                \subsection{The upper bound}

                \begin{Prop1} In the notation of Theorem 1:
                $$\textrm{e}_{\textrm{alg}}(f) \leq \rho(f)$$
                \end{Prop1}

                \paragraph{Proof.}  Let $\mathcal{A} = \{\alpha_{1}, \ldots, \alpha_{n} \}$ be a  generating set for $\pi$.  Let $\mathfrak{R}$ be the $G$-representation variety of $\pi$. Without loss of generality we assume that $G$ is a subgroup of some $GL_{N}$. The set $\mathcal{A}$ determines an embedding $\iota: \mathfrak{R} \to M_{N \times N}^{n}$ given by:
                $$\iota(\rho) = (\rho(\alpha_{1}), \ldots, \rho(\alpha_{n})) $$

                There is an obvious isomorphism $M_{N \times N}^{n} \cong \mathbb{A}^{N^{2}n}$.

                As discussed in the previous section, there is an integer $L$, finitely many functions: $h_{i}: F_{\mathcal{A}} \to K[X_{1}, \ldots, X_{N^{2}r}] \textit{     } (i = 1 \ldots p)$ whose images are all homogeneous of degree at most $L$, and a finite subset $\{w_{1}, \ldots, w_{p}\} \subset \pi$ such that any element of $F[\mathfrak{R}]^{G}$  can be written as a polynomial in $h_{1}(w_{1}), \ldots, h_{m}(w_{m})$.

                Let $w$ be a word of length $l$ in the elements of $\mathcal{A}$. For $1 \leq i \leq p$ we assign to $w$ the function $h_{i}(w)$ on $\mathbb{A}^{N^{2}n}$. If we think of $w$ as an element of $\pi$ and not just a word, we see that this function is an extension of $h_{i}(w)$ from $\mathfrak{R}$ to all of $\mathbb{A}^{N^{2}n}$. Writing out matrix multiplication in coordinates, we see that $h_{i}(w)$ is a homogeneous function of degree at most $lL$.

                Since all of the functions $h_{i}(w)$ are invariant under conjugation, we can deduce that given $w \in \pi$, $|w|_{\mathcal{A},red} = l$, then $h_{i}(w)$ can be written as a homogeneous function of degree at most $lL$ on $\mathbb{A}^{N^{2}n}$.

                Now, given $w \in \pi$, with $|w|_{\mathcal{A},red} = l$, we have the function $h_{i}(w)$ can be written as a polynomial in $h_{1}(w_{1}), \ldots, h_{p}(w_{p})$, each of which is a homogeneous function of degree at least $1$. Since degree is additive under multiplication of homogeneous polynomials, we have that $h_{i}(w)$ can be written as a polynomial of degree at most $lL$ in $h_{1}(w_{1}), \ldots, h_{p}(w_{p})$.

                Define an affine embedding  $\kappa: \mathfrak{X} \to \mathbb{A}^{m}$ by:
                $$\kappa(\chi) = (h_{1}(w_{1})(\chi), \ldots, h_{p}(w_{p})(\chi))$$

                Given an integer $m$, we can write the action of $f^{\circ m}$ in coordinates as:

                $$\kappa \circ f^{\circ m} = (h_{1}(f^{\circ m}(w_{1})), \ldots h_{1}(f^{\circ m}(w_{p})))$$

                From the above discussion, we see that $$\deg(f^{\circ m};\kappa(\mathfrak{X})) \leq L\max(|f^{\circ m}w_{1}|_{\mathcal{A},red}, \ldots, |f^{\circ m}w_{p}|_{\mathcal{A},red})$$

                Therefore, by the definitions of algebraic entropy and spectral radius:

                $$\textrm{e}_{\textrm{alg}}(f) = \limsup_{n \to \infty}\frac{1}{n} \log(\deg(f^{\circ p};\kappa(\mathfrak{X}))) \leq \limsup_{n \to \infty}\frac{1}{n} \log(\max(|f^{\circ m}w_{1}|_{\mathcal{A},red}, \ldots, |f^{\circ m}w_{p}|_{\mathcal{A},red})) \leq \rho(f)  $$

                $\Box$

                \subsection{The lower bound}

                                \begin{Prop2} In the notation of Theorem 1:
                $$\textrm{e}_{\textrm{alg}}(f) \geq \rho(f)$$
                \end{Prop2}

                The proof of the lower bound is more involved than the proof of the upper bound. We begin by recalling some necessary material.

                 \paragraph{Bruhat Tits trees. }Let $K$ be a nonarchimedian complete field of characteristic $0$,
                equipped with a valuation $\nu$. Let $\mathcal{O}_{K}$ be the ring
                of integers, $\mathcal{M}_{K}$ the maximal ideal of
                $\mathcal{O}_{K}$, $k = \frac{\mathcal{O}_{K}}{\mathcal{M}_{K}}$ its
                residue field. Let $q = |k|$ be the number of elements of $k$, and
                let $p$ be its characteristic. For an algebraic group $G$ defined over $K$, let
                $G_{K}$ be the subgroup of $K$ points. $G_{K}$ has a natural action on
                a simplicial complex called the \textit{Bruhat-Tits building of $G_{K}$}. This building
                plays an analogous role to the symmetric space in the archimedian
                case. We will only need to use this theory for $SL_{2}(K)$, in which
                case the building is a regular tree. All of the information that we use can be found in \cite{trees}. Recall that a \textit{lattice}
                in $K^{2}$ is an $\mathcal{O}_{K}$ submodule of the form
                $\mathcal{O}_{K}v \oplus \mathcal{O}_{K}w$, with $v,w \in K^{2}$
                linearly independent. We always denote by $L_{0}$ the so called
                standard lattice: $L_{0} = \mathcal{O}_{K}\left(
                \begin{array}{cc}
                1   \\
                0   \end{array} \right)\oplus \mathcal{O}_{K}\left(
                \begin{array}{cc}
                0   \\
                1   \end{array} \right)$.
                 Lattices $L$ and $L'$ are
                called \textit{homothetic} if $\exists x \in K$ such that $L = xL'$.
                Homothety is an equivalence relation, and we denote the equivalence
                class of the lattice $L$ by $[L]$. We say that two homothety classes
                [L] and [L'] are \textit{incident} if there are representatives
                $L_{1}$, $L_{2}$ of $[L]$ and L' of $[L']$ such that
                $$L_{2} <_{p} L' <_{p} L_{1} $$
                where the symbol $<_{p}$ is read: is a subgroup of index $p$ in. It
                is a simple exercise to check that incidence is a symmetric
                relation.

                 We are now ready to define the Bruhat-Tits building (which we denote by $\mathcal{T}_{K}$) for
                $SL_{2}(K)$. $\mathcal{T}_{K}$ is a graph with a vertex for each
                homothety class of lattices and two vertices connected by an edge if
                the corresponding homothety classes are incident. $SL_{2}(K)$ acts
                on $\mathcal{T}_{K}$ by simplicial automorphisms. We summarize the
                properties of this graph and the $SL_{2}(K)$ action on it that we
                need.

                \begin{bruhattits}
                \begin{enumerate}
                \item $\mathcal{T}_{K}$ is a $\frac{q^{2} -1}{q-1}$ regular tree.
                \item Given $A \in SL_{2}(K)$, its translation length is given by
                $-2\max(\nu(tr(A)),0)$.
                \item The action of $SL_{2}(K)$ is transitive.
                \item $Stab_{SL_{2}(K)}([L_{0}]) = SL_{2}(\mathcal{O}_{K})$,
                $Stab_{SL_{2}(K)}(A[L_{0}]) = A SL_{2}(\mathcal{O}_{K}) A^{-1}$.
                \item The set of connected components of $T_{K} / [L_{0}]$ (i.e. the set of neighbors of
                $[L_{0}]$) can be identified with $\mathcal{P}(k^{2})$, so that the
                action of $SL_{2}(\mathcal{O}_{K})$ on this set of components is
                conjugate to its action on $\mathcal{P}(k^{2})$ (by taking
                conjugates, this statement can be made for each vertex of
                $\mathcal{T}_{K}$).
                \item The axis of a diagonal matrix passes through $[L_{0}]$.

                \end{enumerate}
                \end{bruhattits}

                We now state and prove a technical lemma for bounding algebraic entropy from below. The following two lemmas set up the conditions for using this technical lemma.

                 \begin{lowerbound}
                Let $K$ be a field of characteristic $0$ equipped an absolute value $|.|_{\nu}$,  let $V$ be an affine variety defined over K, let $f$ be an endomorphism of $V$ which is defined over $K$, and let $y \in K[V_{K}]$. If there exists $P_{0} \in V_{k}$ with the following properties:
                \begin{enumerate}
                \item $\exists \epsilon > 0$ such that $ \forall n: |y(f^{\circ n}(P_{0}))|_{\nu} > \epsilon$
                \item $\limsup_{n \to \infty} \frac{1}{n}\log(\log(|y(f^{\circ n}(P_{0}))|_{\nu})) = l$

                \end{enumerate}
                then $l \leq \textrm{e}_{\textrm{alg}}(f)$.
                \end{lowerbound}

                \paragraph{Proof.} Choose an ordered subset $Y = \{y_{1}, \ldots y_{N}\} \subset K[V_{K}]$ such that $Y$ generates $K[V_{K}]$ and $y_{1} = y$. The set $Y$ defines an embedding $V_{K} \hookrightarrow K^{N}$. For the remainder of the proof we will ignore the difference between points in $V_{K}$ and their image under this embedding.

                Given $P \in K^{N}$, let $P^{(i)}$ denote its $i$-th coefficient and $|P|_{\nu} = \max_{i}|P^{(i)}|_{\nu}$. Note that for $P \in V_{K}$, one has that $P^{(1)} = y(P)$.

                 Extend the endomorphism $f$ to an endomorphism of $\mathbb{A}^{N}$ of degree $d$, which we also call $f$. In coordinates we can write $f$ as a vector of polynomials with coefficients in $K$.

                 The function $\frac{|f(P)|_{\nu}}{|P|_{\nu}^{d}}$ is bounded on the set $\{P \in K^{N} : |P|_{\nu} \geq \epsilon \}$. Therefore, $\exists C$ such that:
                 $$|f(P_{0})|_{\nu} \leq C|P_{0}|_{\nu}^{d} $$

                 Suppose first that $d \geq 2$. Iterating $f$ we get:
                 $$|f^{\circ n}(P_{0})|_{\nu} \leq C^{1 + d + \ldots d^{n-1}}|P_{0}|_{\nu}^{d^{n}} = C^{\frac{d^{n} -1}{d-1}}|P_{0}|_{\nu}^{d^{n}} $$
                 Taking logarithms we get:
                 $$\log(|f^{\circ n}(P_{0})|_{\nu}) \leq \frac{d^{n} -1}{d-1}\log(C) + d^{n}\log(|P_{0}|_{\nu}) = d^{n}[\log(|P_{0}|_{\nu}) + \frac{1-\frac{1}{d^{n}}}{d-1} \log(C)] $$
                 Taking logarithms once again, and manipulating further, we get:

                 $$\frac{1}{n} \log(\log(|f^{\circ n}(P_{0})|_{\nu})) \leq log(d) + \frac{1}{n} \log[\log(|P_{0}|_{\nu}) + \frac{1-\frac{1}{d^{n}}}{d-1} \log(C)]$$

 Therefore, $\exists D > 0$ such that:

                 $$\frac{1}{n} \log(\log(|f^{\circ n}(P_{0})|_{\nu})) \leq log(d) + \frac{D}{n}$$

  Now, since $$|y(f^{\circ n}(P_{0}))|_{\nu} = |f^{\circ n}(P_{0})^{(1)}|_{\nu} \leq |f^{\circ n}(P_{0})|_{\nu} $$
                 then:

                 $$l = \limsup_{n \to \infty} \frac{1}{n}\log(\log(|y(f^{\circ n}(P_{0}))|_{\nu})) \leq \limsup_{n \to \infty} log(d) + \frac{D}{n} = \log(d)$$

                Given an integer $q$, we have that:

                $$\limsup_{n \to \infty} \frac{1}{n}\log(\log(|y(f^{\circ qn}(P_{0}))|_{\nu})) = \limsup_{n \to \infty} \frac{q}{qn}\log(\log(|y(f^{\circ qn}(P_{0}))|_{\nu})) = ql$$

                Replacing $f$ by $f^{\circ q}$ in the above discussion, we have that for any $q$ and any extension of $f^{\circ q}$ to an endomorphism of $K^{N}$ (which we also denote by $f^{\circ q})$:
                $$ql \leq \log(\textrm{deg} (f^{\circ q})) $$

                By dividing by $q$ and using the definition of degree relative to a subvariety we have that:
                $$l \leq \frac{1}{q} \log(\textrm{deg}(f^{\circ};V)) $$

                Taking the limits we get:

                $$l \leq \textrm{e}_{\textrm{alg}}(f) $$

                Now assume that $d = 1$. In this case, we get that:
                $$ |f^{\circ n}(P_{0})|_{\nu} \leq C^{n}|P_{0}|_{\nu}$$.

                Taking logarithms twice and dividing by $n$ we get:

                $$\frac{1}{n} \log(\log(|f^{\circ n}(P_{0})|_{\nu})) \leq \frac{\log n}{n} + \frac{1}{n} \log [\log C + \log |P_{0}|_{\nu}] $$
                Taking limits, we get $l \leq 0$. Since algebraic entropy is always non-negative, then we have $l \leq \textrm{e}_{\textrm{alg}}(f)$, as required. $\Box$

                \begin{twoactions}  Let $F_{n} = <\mathcal{S}> = <x_{1} \ldots..
                x_{n}>$ be a free group on $n$ generators. Suppose that $\alpha$ is
                an action of $F_{n}$ on the 2d-regular tree $\mathcal{T}_{2d}$
                satisfying the following conditions:
                \begin{enumerate}
                \item $\alpha(x_{1}), \ldots, \alpha(x_{n})$ are all hyperbolic with
                translation distance $t$.
                \item  There exists a unique vertex $v_{0}$ such that $\{v_{0}\} =
                \mathcal{L}_{i} \bigcap \mathcal{L}_{j}$ for any $i,j$, where
                $\mathcal{L}_{i}$ is the axis of $\alpha(x_{i})$
                \end{enumerate}
                Then given $w \in F_{n}$, the translation length of $\alpha(w)$ is
                $t|w|_{\mathcal{S},red}$.
                \end{twoactions}

                \paragraph{Proof.}

                Notice that since every hyperbolic automorphism with translation
                distance $t$ is a power of an automorphism with the same axis and
                translation length 1, it is enough to prove the lemma for $t=1$.
                Furthermore, by adding hyperbolic automorphisms we can assume $d =
                n$. In this case, we have that $\alpha$ is conjugate to the action
                of $F_{d}$ on the Cayley graph of $F_{d}$ associated to the
                generating set $\{x_{1} \ldots x_{d}\}$, where $v_{0}$ corresponds to the identity element. $\Box$

                \begin{representation}
                There exists a valuation $\nu$ on $\mathbb{Q}$ (resp. $\mathbb{Q}[i]$), and a
                representation $\Psi: \pi \to SL_{2}(\mathbb{Q})$ (resp. $SU_{2}(\mathbb{Q}[i]))$ such that the induced action of
                $\pi$ on $\mathcal{T}_{\nu}$, the Bruhat-Tits tree associated to
                $SL_{2}(\mathbb{Q}_{\nu})$ (resp. $SL_{2}(\mathbb{Q}[i]_{\nu})$) satisfies the conditions of Lemma $5.5$.
                \end{representation}

                \paragraph{Proof.}
                Let $n = 2g$ and let $\pi = <x_{1} \ldots x_{2g}>$. We separate into two cases.
                \paragraph{The $SL_{2}$ case.}
                Let $p$ be any sufficiently large prime (just how large it needs to
                be will be clear from the construction). Let $\nu$ be the $p$-adic
                valuation and let $\mathcal{T}_{p}$ be Bruhat-Tits tree for
                $SL_{2}(\mathbb{Q}_{p})$. Let $[L_{0}]$ be the homothety class of
                the standard lattice. Let $D = \left(
                \begin{array}{cc}
                \frac{1}{p} & 0  \\
                0 & p  \end{array} \right)$. Then $D$ is hyperbolic, and its axis passes
                through $[L_{0}]$. The segment connecting $[L_{0}]$ to $D[L_{0}]$
                passes through the neighbor of $[L_{0}]$ corresponding to  the point
                $\left(
                \begin{array}{cc}
                1   \\
                0   \end{array} \right)$ in $\mathbb{P}^{2}(\mathbb{F}_{p})$. The
                segment connecting $[L_{0}]$ to $D^{-1}[L_{0}]$  passes through the
                neighbor of $[L_{0}]$ corresponding to  the point $\left(
                \begin{array}{cc}
                0   \\
                1   \end{array} \right)$. Let $S$ be any element of
                $SL_{2}(\mathbb{Z})$  for which the set $$\{\left( \begin{array}{cc}
                0   \\
                1   \end{array} \right), S\left( \begin{array}{cc}
                0   \\
                1   \end{array} \right), \ldots S^{n-1}\left( \begin{array}{cc}
                0   \\
                1   \end{array} \right), \left( \begin{array}{cc}
                1   \\
                0   \end{array} \right), S\left( \begin{array}{cc}
                1   \\
                0   \end{array} \right), \ldots S^{n-1}\left( \begin{array}{cc}
                1   \\
                0   \end{array} \right) \}$$ projects to a set of $2n$ different
                points in $\mathbb{Q}\mathbb{P}^{2}$. If $p$ is chosen to be
                sufficiently high, this set will project to $2n$ different points in
                $\mathbb{P}^{2}(\mathbb{F}_{p})$. Note that                $S[L_{0}] = [L_{0}]$.

                 For $i = 1 \ldots n$ define:
                $$\Psi(x_{i}) =S^{(i-1)}  D S^{-(i-1)}$$
                We have that $\Psi(x_{1}) \ldots \Psi(x_{n})$ act hyperbolically on
                $\mathcal{T}_{p}$, they all have the same translation length, and
                the intersection of any two of their axes is precisely $[L_{0}]$.
                The first two assertions follow from the fact that every
                $\Psi(x_{i})$ is conjugate to $D$. The third assertion follows from
                the fact that the axis of $\Psi(x_{i})$ is $S^{i-1}\mathcal{L}$,
                where $\mathcal{L}$ is the axis of $D$, and by construction these
                are $n$ lines that intersect only at $[L_{0}]$.

                \paragraph{The $G = SU(2)$ case.} The construction is almost identical to
                the previous case. We let $K = \mathbb{Q}[i]$. Let $p$ be a
                sufficiently large prime number. There exist integers $a,b$ such
                that $a^{2} + b^{2} = p$. There are two primes of $\mathcal{O}_{K}$
                that lie above $p$, these are $(a + bi), (a-bi)$. Let $\mathfrak{p}
                = (a + bi)$. Let $\nu$ be the $\nu$-adic valuation, and let
                $\mathcal{T}_{\mathfrak{p}}$ be the Bruhat-Tits tree for
                $SL_{2}(K_{\nu})$. Once again, let $[L_{0}]$ be the homothety class
                of the standard lattice.

                Let $D = \left(
                \begin{array}{cc}
                \frac{(a + bi)^{2}}{p} & 0  \\
                0 & \frac{(a - bi)^{2}}{p}  \end{array} \right)$. Then $D$ is hyperbolic
                with axis passing through $[L_{0}]$. If we take $S$ to be any
                element of the $\mathbb{Q}$-points of $SU(2)$ that is not of finite
                order, and whose elements have denominators that are coprime to $p$ then the construction from the previous case may be applied
                verbatim to this case. $\Box$

                \paragraph{Proof of Proposition $5.2$.}

                First the sake of simplicity, we first assume that $G= SL_{2}(\mathbb{R})$, or $SU(2)$.

                Fix $\epsilon \geq 0$. Choose a generating set $\mathcal{S} = \{x_{1}, \ldots, x_{2g} \}$ of $\pi$ for which:

                $$|\limsup_{n\to \infty} \frac{1}{n} \log (|f^{\circ n}x_{1}|_{\mathcal{S},red}) - \rho(f)| \leq \epsilon $$

                Choose a representation $\Psi$ as in Lemma $5.6$, and let $\nu$
                and $\mathcal{T}_{\nu}$ be as in the construction of $\Psi$. Let
                $|.|_{\nu}$ be the absolute value associated to $\nu$, i.e. $|.|_{\nu} = p^{-\nu(.)}$. Let $\psi \in
                \mathfrak{X}$ be the character of $\Psi$. The character $\psi$ is a $\mathbb{Q}$ or $\mathbb{Q}[i]$ point of
                $\mathfrak{X}$.

                By Lemma $5.5$, for any $m$ one has that the translation length
                of $f^{\circ m}(x_{1})$ on $\mathcal{T}_{\nu}$ is equal to
                $|f^{\circ m}(x_{1})|_{\mathcal{S}, red}$. By Proposition $5.3$ part 2 we get:
                $$|f^{\circ m}(x_{1})|_{\mathcal{S}, red} = -2\nu(tr_{f^{\circ m}(x_{1})}(\psi)) $$

                We now wish to apply Lemma $5.4$. In order to set up the notation of the lemma, let $V = \mathfrak{X}$, $K = \mathbb{Q}$ or $K = \mathbb{Q}[i]$ (depending on which part of lemma $5.6$ we used), $|.|_{\nu}$ be the norm defined above, $P_{0} = \psi$, $y = tr_{x_{1}}$, $f = f$.

                Since word length is always positive, we have that $\nu(tr_{f^{\circ m}(x_{1})}(\psi)) < 0$, and thus in the notation of Lemma $5.4$: $|y(f^{\circ m}(P_{0})|_{\nu} \geq 1$.

                Due to our choice of $x_{1}$, we have that:

               $$\limsup_{n \to \infty} \frac{1}{n} \log(\log(|y(f^{\circ n}(P_{0}))|_{\nu})) = \limsup_{n \to \infty} \frac{1}{n} \log(\log(p^{-\nu(tr_{f^{\circ n}x_{1}}(\psi))})) $$

               $$= \limsup_{n \to \infty} \frac{1}{n} \log(\frac{1}{2}|f^{\circ n}(x_{1})|_{\mathcal{S}, red} + p) = \limsup_{n\to \infty} \frac{1}{n} \log(|f^{\circ n}(x_{1})|_{\mathcal{S},red}) \geq  \rho(f) - \epsilon$$

                Thus, by Lemma $5.4$ we have that $\rho(f) - \epsilon \leq \textrm{e}_{\textrm{alg}}(f)$. Since $\epsilon$ was chosen arbitrarily, we get that $\rho(f) \leq \textrm{e}_{\textrm{alg}}(f)$

                Now, suppose $G$ is one of the groups $GL_{N}$, $SL_{N}$, $SU_{N}$, or $SP_{2N}$. Each of these groups contains a copy of $G= SL_{2}(\mathbb{R})$, or $SU(2)$ embedded in the top right corner. If we take the representation $\Psi$ to have image in this copy, and take $y = tr_{x_{1}} - (N-2)$, then the proof proceeds exactly as above.

                If  $G = SO(3)$, notice that $SO(3)$ is double covered by $SU(2)$, and that any $SO(3)$ representation can be lifted to a $SU(2)$ representation where the trace of each element is multiplied by $\pm 1$. Thus, the proof carries over to the $SO(3)$ case. For $G = SO(N)$,  $N \geq 4$ and $G = O(N) (N \geq 3)$ , notice that these groups contain $SO(3)$ embedded as $3 \times 3$ diagonal matrices, and proceed by the same method.  $\Box$

\paragraph{Proof of Theorem 1.} Theorem $1$ is a direct consequence of Propositions $5.1$ and $5.2$. $\Box$

\bibliography{algentropy}
	\bibliographystyle{plain}

\end{document}